\mag=\magstep1
\def\Aut{\operatorname{Aut}}
\def\Sym{\operatorname{Sym}}
\def\Gal{\operatorname{Gal}}
\def\ord{\operatorname{ord}}
\def\val{\operatorname{val}}
\documentstyle{amsppt}
\document
\topmatter
\author
Tomohide Terasoma
\endauthor
\title
Degeneration of curves and analytic deformations
\endtitle
\endtopmatter
 
\heading
\S 0 Introduction
\endheading

Let $\Delta$ be a complex disk $\{z \in \bold C \mid \mid z \mid <1\}$
and $p: \Cal D \to \Delta$ be a proper flat morphism of relative
dimension one.  It is called a minimal degeneration of genus $g$ curves,
if (1) $\Cal D$ is smooth, (2) the restriction $p : p^{-1}(\Delta-\{0\})$
is a smooth family of genus $g$ curves, and (3) the fiber $p^{-1}(0)$
contains no smooth rational curves of the first kind.  By choosing
a metric $\mu$ on $p^{-1}(\Delta-\{ 0\})$, a loop $\gamma$ around $0$ with
a base point $\delta \in \Delta-\{ 0\}$ induces a differentiable 
automorphism $\Gamma (p,\mu ,\gamma )$ of $p^{-1}(\delta)$, whose
class $\Gamma' (p) \in MC(p^{-1}(\delta))$ is independent of the choice
of the metric $\mu$ and the homotopy class of $\gamma$.  Let $S_g$ be
an oriented closed differentiable surface of genus $g$ and 
$D : p^{-1}(\delta ) \to S_g$ be a diffeomorphism.  Then the conjugacy 
class $\Gamma (p)$ of $D\circ \Gamma' (p)\circ D^{-1}$ 
in the mapping class group $MC_g = MC(S_g)$ of 
$S_g$ is independent of the diffeomorphism $D$.
Then $\Gamma (p)$
is known to be a pseudo-periodic with negative Dehn twists.  It is also
known by [MM] that for any pseudo-periodic conjugacy class $\Gamma$
with negative Dehn twists,
there exists a minimal degeneration of  
genus $g$ curves $p: \Cal D \to \Delta$ such that $\Gamma (p) = \Gamma$
unique upto diffeomorphism over $\Delta$.  If $p_1: \Cal D_1 \to \Delta$
are obtained by a smooth holomorphic deformation of 
the degeneration $p_2 : \Cal D_2 \to \Delta$, then $\Cal D_1$ and $\Cal D_2$
are diffeomorphic over $\Delta$ and $\Gamma (p_1)$ is equal to $\Gamma (p_2)$.
In this paper, we prove that if two minimal degeneration of curves
$p_1: \Cal D_1 \to \Delta$ and $p_2:\Cal D_2 \to \Delta$ satisfies 
$\Gamma (p_1)= \Gamma (p_2)$, then these two degenerations 
$p_1, p_1$ are contained
in the same equivalence class generated by smooth holomorphic deformations.

  The contents of this paper is organized as follows.  We introduce general 
notations for stable curves and stable curves with an action of cyclic
group.  Most of the terminology appeared here are introduced in [MM]
in the topological context.  In the next section, we introduce several 
moduli problems related to level structures and group actions.  Here we
use the notion of the type of a group action introduced in \S 1 and
the begining of \S 2.  The main result of this section is the connectedness
of moduli of curves with an action of cyclic group (Proposition 2.2 (2)).  
In \S 3, we study
the group action of the local moduli space.  Here we introduce a coordinate
on local moduli compatible with the action of the cyclic group.
In \S 4, we study the local moduli map attached to the stabilization
of the degeneration of genus $g$ curves.  The screw numbers introduced
by [MM] appears here as an analytic invariant.  The statement and
the main theorem are given in \S 5.

  The author would like to express his thank to S.Takamura and T.Ashikaga 
for discussions.  This work is done during his stay in Max-Planck-Institute
f\"u r Mathematik and the author is appreciated for its hospitality.

\heading
\S 1 A stable curve with an action of a cyclic group
\endheading

In this section, we introduce several notations concerning stable curves with 
an action of a cyclic group.  Let $g$ be a natural number greater 
than 1.  Let $S$ be an analytic space and $f : \Cal C \to S$
be a proper flat morphism of relative dimension one with connected
fibers.  The 
morphism $f$ is called a stable curve of genus $g$ ([DM]), if
\roster
\item
$\dim H^1(f^{-1}(s), \Cal O_{f^{-1}(s)}) = g$
for all $s \in S$.
\item
All the fibers $f^{-1}(s)$ are reduced curves whose singularities are at 
most nodes.  Here a point $p\in f^{-1}(s)$ is called a node if there
exists an open set $U$ of $p$ such that $U$ is isomorphic to
$\{(x,y) \in \bold C^2 \mid \mid x \mid < 1, \mid y \mid <1, xy =0 \}$.
We do not impose that $f^{-1}(s)$ is irreducible for $s \in S$.

The genus of the normalization $\tilde C_i$ of a component $C_i$
in $f^{-1}(s)$ is called the genus of the irreducible component $C_i$.
\item
Any non-singular genus $0$ component meets with another components
at least $3$ points.
\endroster

For any stable curve $C$ of genus $g$ over a point, we can associate
a graph called dual graph $\tau$ of $C$ as follows.  The sets of vertices 
$V(\tau )$ and edges $E(\tau)$
of the graph $\tau$ is given by the set of irreducible components of
$C$ and the singular points of $C$ respectively.  A vertex $v$ is connected
to an edge $e$ if and only if the corresponding singularity $s_e$ is 
an element of the corresponding component $C_v$.
The normalization of the curve $C_v$ is denoted by $\pi_v : \tilde C_v \to
C_v$ and the inverse image of singular points under $\pi_v$ are called
special points of $\tilde C_v$.  The set of special points in $\tilde C$
is denoted by $Sp (C)$.  The genus of $\tilde C_v$ is denoted by $g(v)$.
It is a function from the set of vertices $V(\tau )$ to the set of 
natural numbers.  An edge $e$ has two extreme $e^{(1)}, e^{(2)}$ and 
they are called flags.  The set of flags is denoted by $F(\tau )$
and it is identified with $Sp(C)$.
A flag $f$ is called a tail of $v$ if $f$ is connected to the vertex $v$.
The set of tails of $v$ is denoted by $T(v)$.
The number of the tails of $v$ is denoted by $t(v)$.
A graph is called numbered if the index set of
vertices $V(\tau)$ and flags $F(\tau)$ are numbered by $I_v$ and $I_f$
respectively. The union of $I_v$ and $I_f$ is denoted by $I$.  
\demo{Definition}
Let $\tau$ be a connected graph and $g:V(\tau) \to \bold N$ be a function
from the set of vertices to the set of non-negative integers.
A pair $(\tau , g)$ is called a stable graph if any vertex $v$ with
$g(v)=0$ satisfies $t(v) \geq 3$.
\enddemo

  An automorphism $\sigma$ of $C$ of finite order $m$ induces an 
automorphism of the dual graph $\tau$ of $C$ and it is denoted
by $\sigma_*$.  If we fix a numbering of $\tau$, $\sigma_*$
can be identified with a permutation of the index set $I$.  It is easy
to see that the function $g: V(\tau) \to \bold N$ is preserved by $\sigma_*$.

 Let $G$ be the subgroup of $\Aut (C)$ generated by $\sigma$.  Let
$G(v)$, $G(f)$ and $G(e)$ be the stabilizer of a vertex $v$, 
a flag $f$ and an edge $e$ respectively.  
Then we have $G(e) \supset G(f)$ and $\# [G(e):G(f)]$ is 1 or 2.
An edge $e$ is called amphidrome (resp. non-amphidrome) if
$\#[G(e):G(f)] =2$ (resp. $G(e) = G(f)$).
An element $\sigma \in G(v)$
(resp. $\sigma \in G(f)$ , $\sigma \in G(e)$) acts on the normalization
$\tilde C_v$ of the component $C_v$ (resp. a neibourhood of the 
special point, a neibourhood of the singular
point).
For a point $p$ in the normalization $\tilde C$ of $C$, the character of the 
stablizer $G(p)$ of $p$ induced on the tangent space $T_p\tilde C$
of $\tilde C$ at $p$ is called the local representation at $p$ and
denoted by $\rho (p) : G(p) \to \bold C^{\times}$.
If the point $p$ corresponds to a flag $f$, $\rho_x$ is denoted by
$\rho_f$.
Note that this is an injective homomorphism.  If $p$ and
$q$ are in the same orbit under $G$, we have $G(p) = G(q)$ and
$\rho (p) = \rho (q)$.  
A point $p \in \tilde C$ is called a ramification point
for the action of $G$ if $G(p) \neq 1$.  
The order of the group $G(p)$ is called the ramification index
and denoted by $m_p$.  If $p$ corresponds to a flag $f$,
$m_p$ is written by $m_f$.
The order of $G(e)$ is denoted by $m_e$.  
For the local representation $\rho (p) : G(p) \to 
\bold C^{\times}$ at $p$, put $\rho (\bold e(\frac{1}{m_p}))
=\bold e(\frac{b}{m_p})$. Then since $(b, m_p)=1$, there exists
$0\leq a < m_p$ such that $a \cdot b \equiv 1$ 
$(\text{ mod } m_p )$.
The rational number $\frac{a}{m_e}$ 
is called the valency at $p$ for the action of $\mu_m$ and denoted
by $\val (p)$.
For a flag $f$, the valency of the corresponding special point
is denoted by $\val (f)$ and the restriction of $\val$ to $T(v)$ is
denoted by $\val_v$.
The set of ramification points in 
$\tilde C$ is denoted by $R$. Then $R$ is a finite set.
Denote $R \cap Sp(C) = R_1$ and 
$R-R_1 = R_0$. 
An element of $R_1$ (resp. $R_0$) is called 
a singular ramification point (resp. a smooth ramification).
Let $\tilde C_v$ be a component of the normalization of $\tilde C$.
An element 
$$
r_v = \sum_{q \in (R_0 \cap \tilde C_v)/G(v)}
[\val (q)] \in \oplus_{\alpha \in \bold Q \cap [0,1)}\bold Z [\alpha]
$$
is called the type of smooth ramification for $v$.  The type of the action
$\sigma$ on $C$ is defined by the collection 
$T =(r_v (v\in V(\tau )), \val_v (f) (f \in T(v)))_v$.  It is easy 
to see the following lemma.
\proclaim{Lemma 1.1}
\roster
\item
$$\sum_{q \in (R\cap \tilde C_v)/G(v)} \val (q)$$
is an integer.
\item
$$
\bar g_v =\frac{1}{2}
(\frac{1}{\# G(v)}(2g-2) - \sum_{p \in (R \cap \tilde C_v)/G(v)}
(1-\frac{1}{m_p}))+1
$$
coincides with the genus of the quotient $\tilde C_v/G(v)$, and
it is a non-negative integer.
\endroster
\endproclaim
\demo{Definition}
\roster
\item
For an element 
$r = \sum_{\alpha}u_{\alpha} [\alpha]
\in \oplus_{\alpha \in \bold Q }\bold Z [\alpha ]$, 
$\sum_{\alpha}u_{\alpha}<\alpha >$ is denoted by $< r >$.
\item
For a positive rational number $\alpha$, $Cor(\alpha )= 1-\frac{1}{d}$,
where $d$ is the denominator of $\alpha$.  For an element $r= \sum_{\alpha}
u_{\alpha}[\alpha ]$, we define $Cor (r) = \sum_{\alpha}u_{\alpha}\cdot
Cor(\alpha )$.
\endroster
\enddemo

\heading
\S 2 Several moduli problems, smoothness and connectedness
\endheading

  In this section we study the moduli space of smooth marked curves.
Let $g \geq 0$ and $m \geq 1$ be natural numbers. The module
$\oplus_{\alpha \in \frac{1}{m}\bold Z \cap [0,1)}
\bold Z [\alpha ]$ is denoted by $B(m)$.
Let $r = \sum_{\alpha}u_{\alpha}[\alpha] \in B(m)$ with $u_\alpha \geq 0$,
$T$ be a finite set on which $\mu_m$ acts and $\val : T \to 
\frac{1}{m}\bold Z \cap [0,1)$
be an invariant function under the action of $\mu_m$.

\demo{Definition}
The pair $(r, \val)$ is realized if 
\roster
\item
$<r>+\sum_{\bar f \in T/\mu_m}\val (\bar f)$ is an integer.
\item
$$
\bar g = \frac{1}{2}
(\frac{1}{m}(2g -2)-Cor (r) -\sum_{\bar f \in T/\mu_m}
Cor(\val (\bar f))) +1
$$
is a non-negative integer.
\item
If $\bar g=0$, no proper divisor $m'$ of $m$ has the following 
property:
$$\text{
$t \in B(m')$ and $Im(\val ) \in \frac{1}{m'}\bold Z$
}.
$$
\endroster
\enddemo
From now on, we fix a sufficiently big prime number $l$ such that
(1) $\mu_m \subset \bold F_l^{\times}$, and
(2) $\Aut (C) \to \Aut (Pic^0(C)_l)$ is injective for all curves of genus $g$.
We define a 
symplectic structure and $\mu_m$-action on $V= \bold F_l^{2g}$
as follows.  
The representation of $\mu_m$ on $\bold F_l^{\times}$ defined
by $\mu_m \ni x \mapsto x^t \in \bold F_l^{\times}$ 
($0\leq t < m$) is denoted by $\bold F^l(t)$.
Let $V^{01}$ and $V^{10}$ be $\mu_m$-modules 
defined by 
$V^{01}= \oplus_{t=1}^{m-1}\bold F_l(t)^{h^{01}(t)}$ and
$V^{10}= \oplus_{t=1}^{m-1}\bold F_l(t)^{h^{10}(t)}$, where
$$
\align
h^{01}(t ) =& g-1+<t \cdot r>+\sum_{\bar f \in T/\mu_m}<t \val (\bar f)> \\
h^{10}(t ) =& g-1+<-t \cdot r>+\sum_{\bar f \in T/\mu_m}<-t \val (\bar f)> \\
\endalign
$$
Then there exists a natural perfect pairing
$(,) :V^{01} \times V^{10} \to \bold F_l$.  Let $V_p = V^{01} \oplus V^{10}$.
This vector space has a symplectic structure $\psi_p$ defined by
$$
\align
\psi_p (a,a') = \psi_p (b,b') = 0 
& \qquad (a,a' \in V^{01}, b,b' \in V^{10}) \\
\psi_p (a, b) = -\psi_p (b, a) = (a,b) 
& \qquad (a \in V^{01}, b \in V^{10}). \\
\endalign
$$
Let $\bar V$ be the vector space over $\bold F_l$ with the base
$\alpha_1, \dots, \alpha_{\bar g}, \beta_1, \dots ,\beta_{\bar g}$.
The symplectic form $\bar \psi$ on $\bar V$ is given by
$$
\bar\psi (\alpha_i, \alpha_j) = \bar\psi (\beta_i, \beta_j) = 0,
\bar\psi (\alpha_i, \beta_j )= \delta_{i,j}.
$$
Then the direct sum $V = \bar V \oplus V_p$ has a symplectic form
$\bar\psi \oplus \psi_p$.

Now we define moduli stack $\Cal M_{g,t}$, $\Cal M_{g,t}(l)$, 
$\Cal M_{g,t}(r, \val)$ and $\Cal M_{g,t}(r, \val, l)$.  An object of
$\Cal M_{g,t}$ is a pair $(\Cal C \to S, P_f (f\in T))$ consisting of
a smooth curve $\Cal C \to S$ of genus $g$ over an analytic variety
$S$ and a set of sections $P_f : S \to C$ indexed by the set $T$ of 
the cardinality $t$.
We assume that $P_f(s) \neq P_{f'}(s)$ for all $s \in S$ if $f \neq f'$.
A morphism from $(\Cal C_1\to S_1, P_{1f})$ to $(\Cal C_2\to S_2, P_{2f})$ is 
a following commutative diagram 
$$
\CD
\Cal C_1 @>>> \Cal C_2 \\
@VVV @VVV \\
S_1 @>>> S_2 \\
\endCD
\tag{*}
$$
which is (1) cartesian, and (2) preserving sections indexed by $T$.

An object of $\Cal M_{g,t}(l)$ is a triple 
$$
(\Cal C \to S, P_f (f \in T),  \phi : Pic^0(\Cal C/S)_l \simeq V),
$$
where $\Cal C \to S$ and $P_f$ are as before and $\phi$ is an isomorphis
of local system of $\bold F_l$ vector spaces compatible with the symplectic
structure.  The morphism in $\Cal M_{g,t}(l)$ is defined in the same way.
In this case, we impose that the cartesian product induces an isomorphism 
of local system $Pic^0_l$ compatible with the third data $\phi$.
It is known that $\Cal M_{g,t}(l)$ is connected and
representable if $l\geq 3$.
Let $(r , \val )$ be a realized pair.  An object of
$\Cal M_{g,t}(r, \val )$ is a triple
$(\Cal C \to S, P_f, \psi : \mu_m \to \Aut (\Cal C/S))$ such that
$\Cal C \to S$ and $P_f$ are as before and the action of
$\mu_m$ on $\Cal C$ via $\psi$ preserve that set of sections 
$\{P_f\}_{f \in T}$.
For an action $\psi$, we define valency 
$\val : R \to \frac{1}{m}\bold Z\cap [0,1)$ as before.  We impose that
$r = \sum \val (p)$, and the restriction of $\val$ to $T$ is equal to 
the given map $\val$.  A pair $(r, \val )$ is called the type of the
action $\psi : \mu_m \to \Aut (\Cal C/S)$.
A morphism in $\Cal M_{g,t}(r, \val )$ is a cartesian
product (*) compatible with the action of $\mu_m$.  An object in 
$\Cal M _{g,t}(r, \val, l)$ is a quadraple $(\Cal C, \{ P_f\}
(f \in T), \psi , \phi )$, where $\Cal C \to S$ and $P_f$
$(f \in T)$ are the same as before and 
$\psi : \mu_m \to \Aut (\Cal C)$ and $\phi : Pic^0(\Cal C /S)(l)
\to V$, where the type of the action $\psi$ is $(r, \val )$
and $\phi$ is an isomorphisms of local system compatible with
$\mu_m$-actions and the symplectic structures.

  By forgetting structure, we have the following commutative
diagram:
$$
\CD
\Cal M_{g,t}(r, \val , l) @>{c}>> \Cal M_{g,t}(l) \\
@V{a}VV @VV{b}V \\
\Cal M_{g,t}(r, \val ) @>>{d}> \Cal M_{g,t} \\
\endCD
$$
\proclaim{Proposition 2.1}
\roster
\item
$\Cal M_{g,t}$  and $\Cal M_{g,t}(l)$ are smooth algebraic stack and 
$\Cal M_{g,t}(l)$ is representable.
\item
The morphism $a$ and $b$ are etale.
\item
The morphism $c$ is an immersion. Especially, $\Cal M_{g,t}(r, \val, l)$
is 
\linebreak
representable.
\endroster
\endproclaim
\demo{Proof}
For the statement (1) and the etaleness of the morphism $a$, see [DM]. 
To prove the etaleness of $b$, we consider an object over $S$;
$S \to \Cal M_{g,t}(r, \val )$.  Let 
$(\Cal C \to S, \{ P_f\}_{f \in T}, \phi : \mu_m \to \Aut (\Cal C /S))$
be the corresponding object.  Then $S \times_{\Cal M_{g,t}(r, \val )}
\Cal M_{g,t}(r, \val, l)$ is represented by 
$$
Isom _{S-gp, \text{symplectic},\mu_m\text{-module}}
(Pic^0(\Cal C/S)_l, V),
$$
which is etale over $S$.  This proves the etaleness of $b$.
To prove the locally closedness of $c$, we take an object
$(\Cal C \to S, \{P_f\}, \psi : Pic^0(\Cal C/S)_l \to V)$.
Consider the natural map 
$$
P_l : \Aut_S(\Cal C) \to \Aut_S(Pic^0(\Cal C/S)_l)
$$
They are finite and unramified over $S$.  Therefore $P_l$ is also finite
and unramified.  Moreover it is injective by the choice of $l$.
Therefore, $P_l$ is closed for sufficiently
large $l$.  Therefore the morphism
$$
\align
Hom_S(\mu_m, \Aut_S(\Cal C)) \to &  Hom_S(\mu_m, \Aut_S(Pic^0(\Cal C/S)_l)) \\
\to &  Hom_S(\mu_m, \Aut_S(V)) \\
\endalign
$$
is also closed.  Let $nat$ be the element of $Hom_S(\mu_m, \Aut_S(V))$
corresponding to the action of $\mu_m$ to $V$.  Then 
$Z = (Ad (\phi ) \circ P_l)_*^{-1}(nat)$
is a closed subvariety of $S$.  It is easy to see that 
$\Cal M_{g,t}(r, \val ,l )\times_{\Cal M_{g,t}}S$ is a subfunctor of $Z$.
Consider the corresponding $\mu_m$-action on $\Cal C\mid_Z \to Z$.
Since the type is a constructible function on $Z$, the strata 
$Z(r ,\val )$ corresponding to the type $(r, \val )$ is a
locally closed subvariety of $Z$.  This proves the statement (3).
The stack $\Cal M_{g,t}(r, \val , l)\times_{\Cal M_{g,t}}S$ is represented
by $Z(r, \val )$.  As a consequence, $\Cal M_{g,t} (r, \val)$ is
an algebraic stack.
\enddemo

  Let $USp(V)$ be the group 
$$
\align
\{ \phi \in \Aut (V) \mid &
(\phi v, \phi w) = (v, w) \text{ for all } v,w \in V, \\
& \phi (gv) = g\phi (v) \text{ for all } g \in \mu_m, v \in V\} \\
\endalign
$$
For an object $a=( \Cal C \to S, \{ P_f\},\psi : \mu_m \to \Aut (\Cal C), 
\phi : Pic^0(\Cal C/S)_l \to V)$ in $\Cal M_{g,t}(r, \val, l)$,
and an element $\sigma \in USp(V)$,
we can define an object $\sigma (a)$ in $\Cal M_{g,t}(r, \val , l)$ by
$$
( \Cal C \to S, \{ P_f\},\psi : \mu_m \to \Aut (\Cal C), 
\sigma\circ\phi : Pic^0(\Cal C/S)_l \to V)
$$
This action gives rise to an action of $USp(V)$ on 
$\Cal M_{g, t}(r, \val , l)$.
\proclaim{Proposition 2.2}
\roster
\item
The scheme $\Cal M_{g,t}(r, \val ,l)$ is smooth
\item
The group $USp(v)$ acts transitively to the set of connected
components of $\Cal M_{g,t}(r, \val, l)$.
\endroster
\endproclaim
\demo{Proof}
Let $p$ a point of $\Cal M_{g,t}(r, \val , l)$ and $U$ be 
a sufficiently small neibourhood of $p$.  Let $\Cal C \to U$,
$\phi : \mu_m \to \Aut (\Cal C/U)$, and $\psi : Pic^0(\Cal C/U)_l \to V$
be the corresponding curve $\Cal C$, automorphism of $\Cal C$,
and isomorphism of etale sheaves on $U$.  By taking quotient 
$\Cal D = \Cal C / \mu_m$, we get a genus $\bar g$ curve with 
a level structure $\psi^{\mu_m} : Pic^0(\Cal D/U)_l \simeq
Pic^0(C/U)^{\mu_m} \simeq V^{\mu_m}$.  This gives a morphism
$U \overset{\alpha}\to\to \Cal M_{\bar g}(l)$.
By taking sufficiently small $U$, we may assume that the image 
$\alpha (U)$ is contained in a sufficiently small neibourhood
$W$ of $\alpha (p)$ in $\Cal M_{\bar g}(l)$.  Let $\Cal D_W \to W$ be the 
corresponding curve over $W$. Then we have 
$\Cal D \simeq \Cal D_W \times_W U$.
Let $\Sym (\Cal D_W/W) = \prod \Sym^{u_\alpha}(\Cal D_W/W) \times_W
\Cal D_W^{\bar t}$, where
$r = \sum_{\alpha \in \bold Q}u_\alpha [\alpha ]$ and
$\bar t = \# (T/\mu_m)$.  Then the branch locus of $\Cal D$ gives
a natural morphism $\beta : U \to \Sym (\Cal D_W/W)$ from $U$
to $\Sym (\Cal D_W/W)$.  
More precisely, let $\Cal L_U = (\pi_*\Cal O_{\Cal C})(\omega )$,
where $\omega : \mu_m \to \bold C^{\times}$ is the natural inclusion.
Then $\Cal L_U^{\otimes m} = \Cal O_{\Cal D}(-R)$, where
$$
R = \sum_\alpha (m\alpha )b_\alpha +
\sum_{\bar f \in T /\mu_m}m\val (\bar f)\cdot p_{\bar f},
$$
where $b_\alpha: U \to \Sym^{u_\alpha}(\Cal D_W/W)$ and
$p_{\bar f} : U \to \Cal D_W$.
Let $\tilde W$ be a sufficiently small
neighbourhood of $\beta (p)$ containing the image $\beta (U)$  of $\beta$
by changing $U$ if necessary.  Let $\Cal D_{\tilde W} \to \tilde W$
be the fiber product $\Cal D_W \times_W \tilde W$ and $R$ be the
divisor on $\Cal D_{\tilde W}$ corresponding to the map 
$\tilde W \to \Sym(\Cal D_{\tilde W}/\tilde W )$.  
Since $\tilde W$ is sufficiently small,
we can find a line bundle $\Cal L_{\tilde W}$ on $\Cal D_{\tilde W}$
such that
\roster
\item
$p^*\Cal L_{\tilde W} = ((\pi_p)_* \Cal O_{\Cal C_p})(\chi )$,
where $p : \{ p\} \to U \to \tilde W$ and $\pi_p : \Cal C_p
\to \Cal C_p / \mu_m \simeq \Cal D_p$.
\item
$\Cal L_{\tilde W}^{\otimes m} \simeq \Cal O (- R)$
on $\Cal D_{\tilde W}$.
\endroster
Let $\Cal C_{\tilde W} = Spec( \Cal O_{\tilde W} \oplus \Cal L_{\tilde W}
\oplus \cdots \oplus \Cal L_{\tilde W}^{\otimes (m-1)})$ be the scheme
whose algebra strucure is defined by $\Cal L_{\tilde W}^{\otimes m}
\simeq \Cal O(-R)$.  One can define $\mu_m$-action on $\Cal C_W$
by $\phi (\sigma )^* \mid_{\Cal L_{\tilde W}} = \omega (\sigma)$.
If we take $\tilde W$ to be simply connected, we can find an isomorphism
of etale sheaves $\phi : Pic^0(\Cal C_{\tilde W})_l \simeq V$ and
an isomorphism $p^*\Cal C_{\tilde W} \simeq \Cal C_p$ compatible
with $\phi$.  This defines a morphism $\tilde W \overset{\gamma}\to\to
\Cal M_{g,t}(r, \val , l)$ such that $\gamma\circ\beta : U \to \tilde W
\to \Cal M_{g,t}(r, \val , l)$ is the natural inclusion.  By changing
$\tilde W$ sufficiently small, we may assume $\beta^{-1}(\tilde W) =U$.
In this situation, $\beta$ is an isomorphism and 
$\tilde W$ is smooth.  This proves the smoothness of 
$\Cal M_{g,t}(r, \val,l)$.

(2)  Let $P_1,P_2$ be points in $\Cal M_{g,t}(r, \val, l)$.  We prove
that there exists an element $\sigma$ such that $P_1$ and $\sigma (P_2)$
are connected by a path in $\Cal M_{g,t}(r, \val , l)$.  Let
$(C_i, \{ P_{f i}\}, \phi_i, \psi_i)$ be the quadraple
corresponding to the point $P_i$ ($i=1,2$).  Let $D_i = C_i/\mu_m$
and $E_i$ be the maximal unramified covering of $D_i$.  Let
$\mu_{m'}\subset \mu_m$ be the subgroup of $\mu_m$
corresponding to $E_i$. 
Let $Q_i$ be the point of $\Cal M_{\bar g}(l)$
corresponding to $D_i$ and the isomorphism
$Pic^0(D_i)_l \simeq Pic^0(C_i)_l^{\mu_m}
\simeq V^{\mu_m} = (\bold F_l)^{2\bar g}$.  
Since $\Cal M_{\bar g}(l)$ is connected there exists a path $\gamma_1$
connecting $Q_1$ and $Q_2$.  Taking a lift $\tilde \gamma_1$ of
$\gamma_1$ in $\Cal M_{g,t}(r, \val , l)$, we may assume $Q_1 = Q_2$.
The pair $(E_i, \mu_m/\mu_{m'} \to \Aut ( E_i))$ corresponds to a 
surjective homomorphism from $\pi_1(D_i)$ to $\mu_m/\mu_{m'}$.
Since $T_g \to Sp(2\bar g, \bold F_l) \times
Sp(2\bar g, m\bold Z/m'\bold Z)$ is surjective, we can choose a path
$\gamma_2$ with the base point $Q_1=Q_2$ 
such that the pair $(E_1, \mu_m/\mu_{m'} \to \Aut (E_1))$ 
is analytically continued to the pair 
$(E_2, \mu_m/\mu_{m'} \to \Aut (E_2))$ along the path $\gamma_2$.
By taking a lift $\tilde \gamma_2$ of a path $\gamma_2$ in 
$\Cal M_{g,t} (r, \val , l)$ again, 
we may assume the pair 
$(E_1, \mu_m/\mu_{m'} \to \Aut (E_1))$ is isomorphic to
$(E_2, \mu_m/\mu_{m'} \to \Aut (E_2))$. The curve $E_i$ and $D_i$ 
are denoted by $E$ and $D$ respectively.
Let $f_{i,1}, f_{2}$ and $f_{i,3}$ be the natural morphisms
$$
f_{i,1}: C_i \to E, f_{2}: E \to D, f_{i,3} : C_i \to D.
$$
Let $\omega$ be the natural homomorphism $\mu_m \to \bold C^{\times}$.
Then the character $\omega ' = \omega^{m'}$ is considered as a character
of $\mu_m/\mu_{m'} = \Gal (E/ D)$.  Let 
$\Cal L_i = f_{i,3*}\Cal O_{C_i}(\omega )$ and 
$\Cal A = f_{2*}\Cal O_{E}(\omega ')$.  Then we have the
natural homomorphism
$\Cal L_i ^{\otimes m'} \to f_{i,3*}\Cal O_{C_i}(\omega ')
\to \Cal A$, and this composite morphism is denoted
by $\theta_i$.  Since $f_{2}$ is unramified, $\Cal A^{\otimes (m/m')}
\simeq \Cal O_D$.  Therefore we have a morphism
$\theta_i^{\otimes (m/m')}:\Cal L_i^{\otimes m} \to \Cal O_D$.
This morphism defines an effecitve divisor $\tilde R_i$.
As in the proof of (1), 
$\tilde R_i$ can be written as 
$$
\tilde R_i = \sum_\alpha (m\alpha )b_\alpha +
\sum_{\bar f \in T /\mu_m}m\val (\bar f)\cdot p_{\bar f},
$$
where $b_\alpha \in Sym^{\mu_\alpha}(D)$ and $p_{\bar f} \in D$.
The largest common divisor of $m\alpha$ ($u_\alpha \neq 0$) and
$m \val (\bar f)$ ($ \bar f \in T/\mu_m$) is equal to $m/m'$.
Therefore $\tilde R_i$ is divisible by $(m/m')$ and $(m'/m)\tilde R_i$
is denoted by $R_i$. 
This is equal to the effective divisor
defined by $\theta$ and we have $\Cal A = \Cal L_i^{\otimes{m'}}
\otimes \Cal O_D(R_i)$. 
For an element $\xi = (b_\alpha , p_{\bar f})$, we define $R(\xi )$
by
$$
R(\xi ) = \sum_\alpha (m'\alpha )b_\alpha +
\sum_{\bar f \in T /\mu_m}m'\val (\bar f)\cdot p_{\bar f}.
$$
Let $r = \sum_\alpha m'\alpha u_\alpha + 
\sum_{\bar f \in T/\mu_m}m'\val (\bar f)$.
Let $F_{\Cal A}$ be the fiber product 
\linebreak
$Sym (D) \times_{Pic^{m'r} (D)}
Pic^{r}(D)$ for morphisms 
$Sym(D) \ni \xi \mapsto \Cal O_D(-R(\xi )) \in$ 
\linebreak
$Pic^{m'r}(D)$
and $Pic^{r}(D) \ni \Cal L \mapsto \Cal L^{\otimes m'}\otimes\Cal A^{-1}
\in Pic^{m'r}(D)$.
Then the pair $(R_i, \Cal L_i)$ defines a point of $F_{\Cal A}$.
\enddemo
\proclaim{Proposition 2.3}
The fiber product $F_{\Cal A}$ is 
connected.
\endproclaim
\demo{Proof}
Since $Pic^r (D)$ is an etale covering of $Pic^{m'r}(D)$ associated to 
\linebreak
$m' H_1(Pic^{m'r}(D),\bold Z)$, it is enough to show that the 
composite
$$
\pi_1 (Sym(D)) \to \pi_1(Pic^{m'r}(D)) \to \pi_1(Pic^{m'r}(D)) \otimes
\bold Z /m'\bold Z
\tag{*}
$$
is surjective.  Consider the natural map $D^s \to \Sym (D)$, where
$s =\sum_\alpha u_\alpha + \#(T/\mu_{m})$, and write the composite 
$D^s \to Pic^{m'r}(D)$ as $(d_i)_{i=1, \dots ,s} \mapsto$
\linebreak
$\Cal O_D(-\sum_i\mu_i d_i)$, wehre $\mu_i$ is the multiplicity 
for the $i$-th component.  
To get the 
surjectivity of (*), it is enough to prove the surjectivity of
$$
H_1(D, \bold Z)^{\oplus s} \to H_1(\Sym (D), \bold Z) \to
H_1(Pic^{m'r}(D), \bold Z) \otimes \bold Z /m'\bold Z.
$$
Since the map $H_1(D, \bold Z) \to H_1(Pic^{\mu_i}(D), \bold Z)$
induced by the map $p \mapsto \Cal O(\mu_i p)$ is given by 
$\mu_i$-multiplication, and $gcd(\mu_i,m') =1$,
we get the required surjectivity.  
\enddemo

\heading
\S 3 $\mu_m$-action for a local moduli space.
\endheading

Now we return to the stable curve $C = \cup_{v\in V(\tau)} C_v$.
Consider a $\mu_m$-action on $C$ whose type is $(r_v, \val_v)_v$.
Let $G = \mu_m$ and $G(v)$ be the stabilizer of $v$ in $G$.
Let $X \overset{\pi}\to\to \bar\Cal M_g$ be an etale covering of 
$\bar\Cal M_g$ with a representable stack $X$ and $\bar p: Spec (\bold C)
\to \bar\Cal M_g$ be the point of $\bar\Cal M_g$ correspoiding to $C$.
Then there exists a point $p \to X$ such that $\bar p=\pi\circ p$.
Let $U$ be a neighbourhood of $p \in X$ and $\Cal C \to U$ be the 
corresponding curve on $U$.  
By [DM], the inverse image of $\bar\Cal M_g -\Cal M_g$ under the
map $U \to \bar\Cal M_g$ is a normal crossing divisor $D$ and by taking
sufficiently small $U$, we may assume $D = \cup_{e \in E(\tau )} D_e$.  Let $W$
be a closed analytic set defined by $\cap_{e\in E(\tau)} D_e$.  
It is easy to see
that for a point $w \in W$, the dual graph of the fiber $f^{-1}(w)$
at $w$ is isomorphic to $\tau$.  Moreover, by taking sufficiently small
$U$, the inverse image $\Cal C_W$ of $W$ under the map $f$ has an
irreducible decomposition $\Cal C_W = \cup_{v \in V(\tau)} \Cal C_{W,v}$.
Let $f_{W,v}: \tilde C_{W,v} \to W$ be the normalization of
$\Cal C_{W,v}$ and $\tilde C_{q,v} = f^{-1}_{W,v}(q)$ for $q \in W$.

Let $\bar V(\tau ) \subset V(\tau )$ be a representative of $\mu_m$-orbit
of the set of vertices
$V(\tau )$.  We consider a $\bold F_l$-vector space $V_v$ for each $v \in V$
as in the last section.  For $v \in \bar V(\tau )$, we choose a marking 
$\psi_{p,v} : Pic^0(\tilde C_{p,v})_l \simeq V_v$ compatible with the
action of $G(v)$.  Then the quadraple
$(\tilde C_{p,v}, P_f (f \in T(v)), \phi_v : G(v) \to \Aut (\tilde C_{p,v}),
\psi_{p,v} : Pic^0(\tilde C_{p,v})_l\simeq V_v)$ defines a point 
$p_v \in \Cal M_v (r_v, \val_v, l) \subset \Cal M_v(l)$, where
$\Cal M_v(r_v, \val_v,l) = $
\linebreak
$\Cal M_{g(v),t(v)}(r_v, \val_v, l)$
and $\Cal M_v(l) = \Cal M_{g(v), t(v)}(l)$.
By using isomorphism $\sigma^{i-1}$
$(i=1, \dots ,\#(G/G(v)))$, we define a point 
$\sigma (p_v) \in \Cal M_{\sigma^{(i-1)}(v)}(r_v, \val_v, l)
\subset \Cal M_{\sigma^{(i-1)}(v)}(l)$ and as a consequence,
we define a point 
$$
p' \in \prod_{V(\tau)/G} \prod_{i=1}^{\#(G/G(v))}\Cal M_{\sigma^{i-1}(v)}
(r_v, \val_v, l) \subset \prod_{v \in V(\tau)}\Cal M_v (l).
$$

By choosing a family of marking $\psi_v : Pic^0(\Cal C_{W,v})_l \simeq V_v$
extending $\psi_{v,p}$, we get an etale morphism
$W \to \prod_v \Cal M_v(l)$.  
For a sufficiently small $U$, this map is an open immersion and $p$
is identified with $p'$.
The fiber product 
$W \times_{\prod \Cal M_v(l)} \prod_v\Cal M_v(r_v, \val_v, l)$
is denoted by $Z$.
We can define the action of $G$ on $\prod_v \Cal M_v(l)$ as follows.
An element $\sigma$ in $G(v)$ acts on the space $\Cal M_v(l)$ by 
$(\tilde\Cal C_v , \psi_v) \mapsto (\tilde\Cal C_v, \sigma\circ \psi_v)$.
Therefore on the product, $\Cal M_v \times G$ the group $G(v)$ acts diagonally.
The quotien $(\Cal M_v \times G)/G(v)$ is denoted by 
$Ind_{G(e)}^G\Cal M_v(l)$
and via the action on the second factor, the group $G$ acts on 
$Ind_{G(e)}^G\Cal M_v(l)$.
Using natural isomorphism $\prod_{v \in V(\tau)} \Cal M_v (l) \simeq
\prod_{\bar v \in V(\tau)/G}Ind_{G(v)}^G \Cal M_v(l)$, we get the action of 
$G$ on $\prod_v \Cal M_v(l)$.
The restriction of $f$ to $\Cal C_Z = f^{-1}(Z)$
is denoted by $f_Z$.
$$
\CD
C &\subset& \Cal C_Z &\subset & \Cal C_W & \subset & \Cal  C \\
@VVV @VVV @VVV @ VVV \\
p &\in& Z &\subset& W &\subset &U \\
\endCD
$$

Since $\Cal C \to U$ is a versal deformation of $C$, and
automorphism of stable curves are discrete, for an 
automorphism $g \in \Aut (C)$,  there exists a sufficiently small
neighbour hood $U_1$ of $p$ and a map $g_* : U_1 \to U$ such that
the automorphism $g$ is induced by the pull of $\Cal C \to U$
by $g_*$.  Since $\Aut (C)$ is a finite group, we may assume that
$U_1$ is stable under the map $g_*$ , $g \in \Aut (C)$.
This defines an action of $\Aut (C)$ on $U$ and by restricting this
representation,
we get an action of $G$ on $U$.
It is easy to see that the subspace $W$ is stable under the action
of $G$.  Moreover, the restriction of this action to $W$
coincides with the action on $\prod_v\Cal M_v(l)$
given before.
Therefore, the subspace $Z$ is fixed part of $W$ under the action of 
$G = <\sigma>$.

We introduce a coordinate of $U$ as follows.
Since $\sigma (D_e) = D_{\sigma (e)}$, and $\cup_{e \in E(\tau )}D_e$
is normal crossing, we can take local equotions $t_e$ of $D_e$
such that $t_{\sigma (e)}$ is a constant multiple of $\sigma^*(t_e)$.
By this choice of  $\{t_e\}_{e \in E(\tau)}$, 
$$
G(e) \ni \sigma \mapsto \frac{\sigma^*(t_e)}{t_e} \in \bold C^{\times}
$$
defines a character $\chi_e$ of $G(e)$.  Since the subspace
$W$ is smooth, we can take a coordinate $\{ t_i\}_{i = 1, \dots , \dim W}$
such that $\sigma^*(t_i) = \chi_i(\sigma )t_i$ , $(\sigma \in G)$
for some character $\chi_i$ of $G$.  As a whole, $\{ t_i\}_{i= 1, \dots , 
\dim W} \cup \{ t_e\}_{e \in E(\tau )}$ forms a local coordinate for
$U$.  The character $\chi_e$ is computed from $(r_v, val_v)_v$ as
follows.

  Let $e^{(1)}$ and $e^{(2)}$ be the two extremes of the edge $e$ and
$u_1$ and $u_2$ be the local coordinate of $\tilde C$ corresponding
to $e^{(1)}$ and $e^{(2)}$.  We may assume that for any 
$\sigma \in G(e^{(i)})$, $\sigma^*(u_i)$ is a constant multiple of
$u_i$.  For a small neighbourhood of $p_e \in \Cal C$, the local equation
is written as $u_1u_2 = t_e$. The group $G(e^{(1)}) = G(e^{(2)})$
is denoted by $G(f)$.

(1) Amphidrome case.  Suppose that there exists $\sigma \in G(e)$
such that $\sigma (e^{(1)}) = e^{(2)}$.  In this case, $\# [G(e) : G(f)] =2$
and for $\sigma \in G(f)$, $\chi_f (\sigma ) = \sigma ^*(u_1)/u_1
= \sigma ^*(u_2)/u_2$ defines a character of $G(f)$. For $\sigma \in G(e)$,
we have
$$
\chi_e(\sigma) = \frac{\sigma^*(t_e)}{t_e} = \frac{\sigma^*(u_1u_2)}{u_1u_2}
= \rho_f(\sigma^2)
$$

(2) Non-amphidrome case.  If $G(e)= G(f)$, $\chi_{e^{(1)}}(\sigma )
= \sigma^*(u_1)/u_1$ and $\chi_{e^{(2)}}(\sigma )
= \sigma^*(u_2)/u_2$ defines a characters of $G(e)$.  By the same argument,
we have $\chi_e = \rho_{e^{(1)}} \rho_{e^{(2)}}$.

\heading 
\S 4 Period of degeneration and monodromy
\endheading

  In this section, we study the period map for the stabilization of a
degeneration of curves.
Let $p :\Cal D \to \Delta$ be a proper morphism of dimension one
whose restriction to $p^{-1}(\Delta-\{0\})$ is smooth. 
If the fiber $p^{-1}(0)$ at 0 contains no rational curves of the first kind,
it is called a minimal degeneration of genus $g$ curves.
Then by the stable reduction
theorem [DM], there exists a covering 
$\pi_m :\Delta_m \to \Delta; t \mapsto t^m =\tau$ 
of degree $m$ such
that the fiber product $\Cal D_m = \Cal D \times_{\Delta}\Delta_m$ 
has stable reduction.
Let $\Cal C \to \Delta_m$ be the stabel model of $\Cal D_m$.  The special
fiber $C$ of $\Cal C$ defines a point $p$ of $\bar\Cal M_g$.  Let $U$ be
a neighbourhood of $p$ which is representable and sufficiently small.
The corresponding curve on $U$ is denoted by $\Cal C_U \to U$.
By the functoriality of stable model, the action of $\Gal (\Delta_m /\Delta)$
acts on $\Cal C$ which is compatible with the natural action on $\Delta_m$.
Restricting this action to the special fiber, we get an action of 
$\Gal (\Delta_m/\Delta)$ on $C$.  By the versality of $\Cal C_U \to U$,
the group $\Gal (\Delta_m /\Delta )$ acts on $U$ by taking sufficiently
small $U$ as in Section 3.  We get the following cartesian diagram
compatible with the action of $\Gal (\Delta_m/ \Delta)$ by changing
$U$ by sufficiently small neighbourhood of $p$.
$$
\CD
\Cal C @>>> \Cal C_U \\
@V{p_m}VV @VV{p_U}V \\
\Delta_m @>>{\Phi}> U \\
\endCD
$$
We choose a coordinate $(t_i, t_e)$ introduced in Section 3.
Using this coordinate the morphism $\Phi$ can be written as
$\Phi (t) = (f_i(t), f_e(t))$, where $f_i$ and $f_e$ are holomorphic
function on $\Delta_m$.  Since the morphism $\Phi$ is equivariant
under the action of $G = \Gal(\Delta_m/\Delta )$, we have the following 
functional equation for $f_i$ and $f_e$.
$$
f_i(\sigma t) = \chi_i(\sigma )f(t),
f_e(\sigma t) = c_{\sigma ,e}f_{\sigma e}(t), 
$$
where $c_{\sigma ,e}= \sigma^*{t_e}/t_{\sigma e}\in \bold C^{\times}$.
As a consequence, there exist germs of holomorphic functions
$\tilde f_i$ and $\tilde f_e$ such that
$$
f_i (t) = t^{e_i}\tilde f_i(t^m),
f_e(t)  = t^{m_e b_e}\tilde f_e(t^{m_e}), 
$$
where $0\leq b_e<1$ and $\rho_f(\bold e(\frac{1}{m_f}))= 
\bold e (b_e)$ if $e$ is amphidrome and
\linebreak 
$\rho_{e^{(1)}}(\bold e(\frac{1}{m_{e^{(1)}}}))
\rho_{e^{(2)}}(\bold e(\frac{1}{m_{e^{(2)}}})) =
\bold e (b_e)$ if $e$ is non-amphidrome.
  Now we consider the minimal resolution of $\Cal C /G$.
Note that the singularity of $\Cal C$ is contained in singular
locus of the special fiber $C$.  Let $x$ be a point in $\Cal C$
fixed by some non-trivial element of $G$, i.e. $G(x) = Stab_x(G) \neq 1$.

(1) The case where $x$ is contained in the smooth part of $C$.
In this case, $G(x)$ acts on the tangent space of $C$ at $x$ via the
local representation $\rho_x$.  Therefore the action of $G(x)$
on the tangent spece of $\Cal C$ is equivalent to the direct sum
$\rho_x \oplus nat$ of $\rho_x$ and the natural representation $nat$.

Therefore the resolution process depends only on $\val (x)$ and
does not depend on the map $\Phi$.

(2) The case where $x$ is contained in the singular locus of $C$.
The local equation of $\Cal C$ at $x$ is 
given by
$$
u_1u_2 = t^{m_eb_e}\tilde f_e(t^{m_e})
$$
in the spce $(u_1, u_2, t)$, 
The action of $G(x) = G(e)$ is given as follows.

(2-1) Amphidrome case.  
Let $\sigma =\bold e (\frac{1}{m_e})$ be the generator 
of $G(e)$ (see \S 3).  Then by changing coordinate,
the action of $\sigma$ is given as 
$$
u_1 \mapsto u_2, u_2 \mapsto \chi_e(\bold e(\frac{1}{m_e}))u_1,
t \mapsto \bold e(\frac{1}{m_e})t.
$$
Therefore the resolution process depends only on $\chi_e$ and the order
of $\tilde f_e$ with respect to the parameter $\tau_e = t^{m_e}$.

(2-1) Non-amphidrome case.  Using the same notation as in (2-1), 
the action of $G_x$  is given as
$$
u_1 \to \chi_{e^{(1)}}(\sigma ) u_1,
u_2 \to \chi_{e^{(2)}}(\sigma ) u_2,
t \mapsto \bold e(\frac{1}{m_e})t.
$$
In this case, the resolution process also depends only on $\chi_{e^{(1)}}$,
$\chi_{e^{(2)}}$ and the order of $\tilde f_e$ with respect to
the parameter $\tau = t^{m_e}$.
\proclaim{Proposition 4.1}
Let $\Phi_j : \Delta_m \to U$ $(j=1,2)$ be a holomorphic map
which satisfies $\Phi_j(\sigma t) = \sigma^*(\Phi (t))$ and
write $\Phi_i = (f_i^{(j)}, f_e^{(j)})$ be using the coordinate
$(t_i, t_e)$.  Suppose that $\ord_t f_e^{(1)}(t) =
\ord_t f_e^{(2)}(t)$ for all $e \in E(\tau)$.  Then there exists a smooth
family $h : \tilde\Cal D \to \Delta$ and $\pm \epsilon \in \Delta$
such that $\Cal D_1 = h^{-1}(-\epsilon )$ and 
$\Cal D_2 = h^{-1}(\epsilon )$ are the minimal resolution of 
$\Phi_1^*(C_U)/G$ and $\Phi_2^*(C_U)/G$ respectively.
\endproclaim
\demo{Proof}
By the assumption, $f_i^{(j)}(t)$ and $f^{(j)}_e(t)$ can be written
as
$$
f_i^{(j)} (t) = t^{e_i}\tilde f_i^{(j)}(t^m),
f_e^{(j)}(t)  = t^{m_e b_e}\tilde f_e^{(j)}(t^{m_e}).
$$
By the assumption we can take two variable function
$\tilde F_i(\tau ,u)$ and $\tilde F_e(\tau,u)$ on 
$\Delta \times \Delta$ such that
$$
\align
\tilde F_i (t, -\epsilon ) = \tilde f_i^{(1)}(t),
&\tilde F_i (t, \epsilon ) = \tilde f_i^{(2)}(t),  \\
\tilde F_e (t, -\epsilon ) = \tilde f_e^{(1)}(t),
&\tilde F_e (t, \epsilon ) = \tilde f_e^{(2)}(t), \\
\endalign
$$
and the order of $\tilde F_e(\tau ,u)$ with respect to $\tau$ is constant
for all $u \in U$ and if $e$ and $e'$ are in the same orbit
under the action of $\mu_m$, $F_e$ is a constant multiple of $F_{e'}$.
By pulling back by the morphism 
$$
\Delta_m \times \Delta \ni
(t,u) \mapsto (t^{e_i}\tilde F_i(t^m,u), t^{m_eb_e}\tilde F_e(t^{m_e},u)) 
\in U,
$$
we get a family of curves on $\Delta \times \Delta$ by taking quatient
and take a resolution, we get the required smooth family 
$\tilde D \to \Delta$.
\enddemo

\heading 
\S 5  Mapping class group and the main theorem
\endheading
 
First we recall several definitions of mapping class group.
Let $\Cal D \to \Delta$ be a minimal degeneration of genus $g$
curves.  Let $0< \epsilon < \delta <1$.  By choosing a metric on
$p^{-1}(\Delta - \Delta_{\epsilon})$, and a path around 0 with the
base point $\delta$, we obtain a $C^{\infty}$ automorphism
of $p^{-1}(\delta )$ and it defines an element 
$\Gamma' (p) \in MC(p^{-1}(\delta ))$
of the mapping class group of $p^{-1}(\delta )$.  Let $S_g$ be 
an oriented closed $C^{\infty}$ surface.  By taking a diffeomorphism
$D : p^{-1}(\delta) \to S_g$ from $p^{-1}(\delta )$ to $S_g$,
we get an element $\Gamma_D (p) = D \circ \Gamma (p) \circ D^{-1}$ 
of the mapping class group $MC_g = MC(S_g)$ of $S_g$.  
The conjugacy class of $\Gamma_D(p)$ does not depend on the choice
of $\epsilon, \delta$ and $D$ and it is denoted by $\Gamma (p)$.

  We review several results of [MM] from the analytic point
of view.  Let $\Cal D \to \Delta$ be a minimal degeneration
of genus $g$ curves.
Let $m$ be the minimal degree of $\Delta_m \to \Delta$
for which $\Cal D \times_{\Delta}\Delta_m$ has stable reduction.
By stable reduction theorem, it is equal to the minimal $m_i$
such that $\det (1- x \Gamma' (p)^m\mid H_1(p^{-1}(\delta ), \bold Z))
= (1-x)^{2g}$.  Therefore $m$ depends only on the conjugacy class 
$\Gamma (p)$ of $\Gamma' (p)$ in the mapping class group.
Let $p_m :\Cal C\to \Delta_m$ be the stable model of 
$\Cal D \times_{\Delta}\Delta_m$
and $C$ be the special fiber $p_m^{-1}(0)$ of $\Cal C$.  Then by the 
functoriality of the stable model the action of $G= \Gal (\Delta_m /\Delta )$
extends to the action of $\Cal C$ and as a consequence, we have an
action of $G$ on the closed fiber $C$ and the smooth part $C^0$ of
$C$.  It is easy to see that $C^0$ is homeomorphic to
$p^{-1}(\delta ) - Cir$, where $Cir$ is the minimal set of simple closed
curve on $p^{-1}$ such that the restriction of $\Gamma (p)$ to
$p^{-1}(\delta ) - Cir$ is periodic.  
Moreover under this homeomorphism,
the action of $G$ on $C^0$ is homotopically equivalent to
that of $\Gamma (p)$ on $p^{-1}(\delta )-Cir$. Therefore the valency
defined in \S 2 is equal to that given in [MM].
Let $U$ be the versal deformation of the special fiber $C$.
We use the same notation $(t_i, t_e)$, $(f_i(t),f_e(t))$, etc. as
in \S 3.
Let $\ord_t(f_e(t))$ be the order of $f_e$ with respect to the parameter
$t$,  Then $\ord_t (f_e(t))/m_e$ of $f_e$ coincides 
with the screw number defined in [MM] and depends only on the conjugacy
class of the mapping class group $\Phi \in \Aut (\pi_1 (\Cal D_{\eta}))$
arising from the family $\Cal D \to \Delta$ of curves.

Now we can prove the following main theorem.
\proclaim{Theorem 5.1}
Let $p_i : \Cal D_i \to \Delta$ ($i=1,2$) be degenerations of genus
$g$ curves.  If $\Gamma (p_1) = \Gamma (p_2)$, then there exists
a sequence of 
proper flat morphisms $e_i : \Cal E_i \to \Delta \times D_i$
($i = 1, \dots , k$) of dimension 1 such that
\roster
\item
The composite $pr_2\circ \Cal C_i : \Cal E_i \to D_i$ is smooth.
\item
The restriction $e_i^{-1}(\Delta^0 \times D_i) \to
\Delta^0 \times D_i$ is smooth.
\item
There exists an isomorphism
$(pr_2\circ e_i )^{-1}(\epsilon ) \overset{\simeq}\to\to
(pr_2\circ e_{i+1} )^{-1}(-\epsilon )$
compatible with the projection to $\Delta$.
\item
$(pr_2\circ e_1)^{-1}(-\epsilon ) \to \Delta$ and 
$(pr_2\circ e_i)^{-1}(\epsilon ) \to \Delta$ are isomorphic to
$p_1 :\Cal D_1 \to \Delta$ and $p_2: \Cal D_2 \to \Delta$
respectively.
\endroster
In other words, if two degenerations of curves are topologically 
isomorphic ot each other, they are equivalent under analytic deformations.
\endproclaim
\demo{Proof}
By the assumption the minimal degree of $\Delta_m \to \Delta$
for which $\Cal D_1 \times_{\Delta}\Delta_m$ has stable reduction
and that for
$\Cal D_2 \times_{\Delta}\Delta_m$ coincides and we denote it $m$.
Let $\Cal C_i$ be the stable model of $\Cal D_i \times_{\Delta}\Delta_m$
and $C_i$ be the special fiber of $\Cal C_i$.  Then the stable graph
$(\tau_i, g)$ associated to $C_i$ are isomorphic to each other.
Moreover the type $(r_{i,v}, \val_{i,v})_v$ of the ramificaiton 
for $C_i$ for the action of $\mu_m$ is also equal and it is denoted 
by $(r_v, \val_v)_v$.  The point in $\bar\Cal M_g$ defined by the special
fiber $C_1$ and $C_2$ are denoted by $p_1$ and $p_2$ respectively.
By Proposition 2.2, there exist liftings $\tilde p_i$ of $p_i$
which belong to the same connected component $K$ of 
$\prod_{v \in V(\tau)/\mu_m}\Cal M_v(r_v, \val_v, l)\subset \bar\Cal M_g(l)$. 
Let $U_i$ ($i=1, \dots, k$) be a sequence of
open sets of $\bar\Cal M_g(l)$ with
$U_i \cap U_{i+1}\cap K \neq \emptyset$ for $i= 1, \dots ,k-1$ and
$\tilde p_1 \in U_1\cap K$ and $\tilde p_2 \in U_k\cap K$.
Since the screw number for $\Cal C_1 \to \Delta_m$ and
$\Cal C_2 \to \Delta_m$ coincides for all $e \in E(\tau)$,
we get the required sequence of morphism $e_i : \Cal E_i \to \Delta
\times D_i$ by Proposition 4.1.
\enddemo

\enddocument